\documentclass[12pt]{amsart}
\usepackage{amssymb}
\usepackage[dvips]{graphicx}
\usepackage{bbm}

\numberwithin{equation}{section}

\newtheorem{thm}{Theorem}
\newtheorem{prop}[thm]{Proposition}
\newtheorem{lemma}[thm]{Lemma}
\newtheorem{cor}[thm]{Corollary}

\theoremstyle{definition}

\newtheorem{example}[thm]{Example}
\newtheorem{remark}[thm]{Remark}

\newtheorem{definition}[thm]{Definition}
\newtheorem{conjecture}[thm]{Conjecture}

\numberwithin{thm}{section}

\newcommand{\card}[1]{{\mid\! #1 \!\mid}}

\newcommand{\N}{\mathbb{N}}
\newcommand{\calC}{\mathcal{C}}

\newcommand{\Z}{\mathbb{Z}}

\newcommand{\R}{\mathbb{R}}

\newcommand{\V}{\mathcal{V}}
\newcommand{\height}{{\rm ht}}
\def\MR{M_\R}

\DeclareMathOperator{\codeg}{codeg}

\DeclareMathOperator{\supp}{supp}

\DeclareMathOperator{\Vol}{Vol}
\DeclareMathOperator{\conv}{conv}
\DeclareMathOperator{\aff}{aff}

\title[$h^*$-polynomials with given degree and linear coefficient]
{Lattice polytopes having $h^*$-polynomials with given degree and linear coefficient}

\author{Benjamin Nill}

\address{Research Group Lattice Polytopes, FU Berlin, Germany}
\email{nill@math.fu-berlin.de}

\begin{document}

\begin{abstract}
The $h^*$-polynomial of a lattice polytope is the num\-erator of the generating function of the Ehrhart polynomial. 
Let $P$ be a lattice polytope with $h^*$-polynomial of degree $d$ and with linear coefficient $h^*_1$. 
We show that $P$ has to be a lattice pyramid over a lower-dimensional lattice polytope, if the dimension of $P$ is greater or equal to 
$h^*_1 (2d+1) + 4d-1$. 
This result has a purely combinatorial proof and generalizes a recent theorem of Batyrev. 
As an application we deduce from an inequality due to Stanley
that the volume of a lattice polytope is bounded by a function depending only
on the degree and the two heighest non-zero coefficients of the $h^*$-polynomial.
\end{abstract}

\maketitle

\section{Introduction and main results}

Let $M$ be a lattice, and $P \subseteq \MR = M \otimes_\Z \R$ be an $n$-dimensional lattice polytope, i.e., the set of vertices of $P$, here denoted by $\V(P)$, 
is contained in the lattice $M$. Throughout, the normalized volume $\Vol(P)$ with respect to $M$ is refered to as the volume of $P$. 
Moreover, two lattice polytopes $P \subseteq \MR$ and $P' \subseteq M'_\R$ are called isomorphic, 
if there is an affine lattice isomorphism $M \cong M'$ mapping $\V(P)$ onto $\V(P')$.

\medskip

Due to Ehrhart and Stanley \cite{Ehr77,Sta80,Sta86} the generating function enumerating the number of lattice points in 
multiples of $P$ is a rational function of the following form:
\[\sum\limits_{k \geq 0} \card{(k \Delta) \cap M} \, t^k =  \frac{h^*_0 + h^*_1 t + \cdots + h^*_n t^n}{(1-t)^{n+1}},\]
where $h^*_0, \ldots, h^*_n$ are non-negative integers 
satisfying the conditions $h^*_0 =1$, $h^*_1 = \card{P \cap M} - n - 1$ and $h^*_0 + \cdots + h^*_n =\Vol(P)$.

\begin{definition}{\rm 
The polynomial $h^*_P(t) := h^*_0 + h^*_1 t + \cdots + h^*_n t^n$ is called the {\em $h^*$-polynomial} of $P$ (see \cite{Bat07,BN07,Sta93}) 
or $\delta$-{\em polynomial} (see \cite{Hib94}). 
The degree of $h^*_P(t)$, i.e., the maximal $i \in \{0, \ldots, n\}$ with $h^*_i\not= 0$, is called the {\em degree} $\deg(P)$ of $P$. 
We define the {\em codegree} of $P$ as $\codeg(P) := n+1-\deg(P)$. 
}
\end{definition}

The geometric meaning of the codegree, introduced by Batyrev in \cite{Bat07}, is given by the following observation:
\[\codeg(P) = \min(k \geq 1 \,:\, k P \text{ has interior lattice points}).\]

The notion of the degree of a lattice polytope was defined in \cite{BN07}, where it was noted that $\deg(P)$ 
should be considered as the "lattice dimension" of $P$. 
This interpretation of the degree was motivated by the following three basic properties: 
First, $\deg(P) = 0$ if and only if $\Vol(P) = 1$. So the unimodular simplex is the only lattice polytope with degree zero. Second, 
by Stanley's monotoni\-city theorem \cite{Sta93} it holds 
$h^*_Q(t) \leq h^*_P(t)$ coefficientwise for lattice polytopes $Q \subseteq P$. 
In particular this implies that the degree is monotone with respect to inclusion. 
For the third property let us recall the notion of lattice pyramids \cite{Bat07}:

\begin{definition}{\rm Let $B \subseteq \R^k$ be a lattice polytope with respect to $\Z^k$. 
Then $\conv(0, B \times \{1\}) \subseteq \R^{k+1}$ is a lattice polytope with respect to $\Z^{k+1}$, 
called the ($1$-fold) standard pyramid over $B$. Recursively, 
we define for $l \in \N_{\geq 1}$ in this way the 
$l$-fold standard pyramid over $B$. As a convention, the $0$-fold standard pyramid over $B$ is $B$ itself.

Let $P,Q \subseteq \MR$ be lattice polytopes with $Q \subseteq P$. 
We say $P$ is a {\em lattice pyramid} over $Q$, if $P \subseteq \MR$ is isomorphic to the $(\dim(P)-\dim(Q))$-fold standard pyramid 
over a lattice polytope $B$, where this isomorphism maps $Q$ onto $B$.
}
\end{definition}

Now, for lattice polytopes $Q \subseteq P$ we observe that $P$ is a lattice pyramid over $Q$ if and only if $\Vol(P) = \Vol(Q)$, or equivalently, 
$h^*_P(t) = h^*_{Q}(t)$. This implies as a third property the invariance of the degree under lattice pyramid constructions. 

\medskip

In \cite{Bat07} Batyrev showed the following theorem:

\begin{thm}[Batyrev]
Let $P \subseteq \MR$ be an $n$-dimensional lattice polytope of volume $V$ and degree $d$. 
If 
\[n \geq 4d {2d+ V -1 \choose 2d},\]
then $P$ is a lattice pyramid over an $(n-1)$-dimensional lattice polytope.
\label{bat-theo}
\end{thm}

Recursively, we see that any lattice polytope $P$ is a lattice pyramid over a lattice polytope $Q$ with $h^*_P(t) = h^*_Q(t)$, where 
the dimension of $Q$ is bounded by a function depending only on the degree and the volume of $P$. Since by \cite{LZ91} there is up to isomorphisms only a finite number of $n$-dimensional 
lattice polytopes with volume $V$, if $n$ and $V$ is fixed, we get the following corollary:

\begin{cor}[Batyrev]
There is only a finite number of lattice polytopes of fixed degree $d$ and fixed volume $V$ up to isomorphisms and lattice pyramid constructions.
\label{bat-cor}
\end{cor}

\medskip

Here, we improve the bound in Batyrev's theorem to the presumably correct asymptotic behaviour:

\begin{prop}
Let $P \subseteq \MR$ as in Theorem \ref{bat-theo}. If 
\[n \geq (V-1) (2d+1),\]
then $P$ is a lattice pyramid over an $(n-1)$-dimensional lattice polytope.
\label{vol-prop}
\end{prop}

Note that for $d=1$ this yields the assumption $n \geq 3 (V - 1)$, while Batyrev's theorem needs $n \geq 2 (V+1) V$.
Since all lattice polytopes of degree one were classified in \cite{BN07}, it could be observed in \cite[Prop. 4.1]{Bat07} that
$n \geq V+1$ is the optimal bound.

\begin{example}{\rm
Here is an example of a lattice polytope with degree $d \geq 2$, volume $V=2$, and dimension $n=2d-1$ that is not a
lattice pyramid: the simplex with vertices $e_0 - e_n, e_1 - e_n, \ldots, e_{n-1} - e_n, e_0 + \cdots + e_{n-1} + (3-2d) e_n$,
where $e_0, \ldots, e_n$ is a lattice basis of $\Z^{n+1}$. The $h^*$-polynomial equals $1 + t^d$. Though this example does not show
that the bound given in Proposition \ref{vol-prop} is sharp, we see again that the asymptotics seems to have the right order.
}
\end{example}

While Batyrev's proof involved commutative and homological algebra, 
our methods are elementary and purely combinatorial.

\bigskip

The main result of this paper shows that the qualitative statement of Theorem \ref{bat-theo} still holds, 
when we fix instead of the volume of $P$ only the "relative" number of vertices $\card{\V(P)} - n-1$, which is an invariant depending only on the 
combinatorics of $P$:

\begin{thm}
Let $c,d \in \N$. Let $P \subseteq \MR$ be an $n$-dimensional lattice polytope having $\leq c+n+1$ vertices and degree $\leq d$. 
If 
\[n \geq  c (2d+1) + 4d-1,\]
then $P$ is a lattice pyramid over an $(n-1)$-dimensional lattice polytope.
\label{main-theo}
\end{thm}

Note that for $d=1$ this yields the assumption $n \geq 3 (c + 1)$, while the optimal bound is
$n \geq 3$ for $c=0$ and $n \geq c+2$ for $c > 0$ by the classification \cite{BN07}.

Now, since $\card{\V(P)} - n - 1 \leq \card{P \cap M} - n - 1  = h^*_1$, we see that the implication in Theorem \ref{main-theo} holds
for $c=h^*_1$. This result motivates the following more general conjecture:

\begin{conjecture}{\rm 
Let $c$, $d$, and $i \in \{1, \ldots, d\}$ be fixed. Then there is a function $f_i$ depending only on $c$ and $d$ 
such that any $n$-dimensional lattice polytope $P$ with $h^*_i=c$ and degree $\deg(P)=d$ is a lattice pyramid over an $(n-1)$-dimensional lattice polytope, 
if $n \geq f_i(c,d)$.
\label{neu}
}
\end{conjecture}

\begin{remark}{\rm 
For $i=1$, the conjecture holds, as we have just seen. 
For $i=2, \ldots, d-1$, the conjecture would follow from the inequalities $h^*_1 \leq h^*_i$. 
In the case of $d=n$, these were proven by Hibi \cite{Hib94}. The author is not aware of any counterexamples for arbitrary degree.
\footnote{In the meantime Henk and Tagami provided a counterexample \cite{HT07}.}

For $i=d$, this conjecture is equivalent to Conjecture 4.2 in \cite{Bat07}, saying that $\Vol(P)$ should be bounded by a function 
in $d$ and $h^*_d$. To see this equivalence we note that $h^*_d > 0$ equals the number of interior lattice points in $\codeg(P) P$, and 
due to Hensley \cite{Hen83} 
the volume of any $n$-dimensional lattice polytope with $l > 0$ interior lattice points is bounded by a function depending only on $n$ and $l$. 
Actually, by a result of Lagarias and Ziegler \cite{LZ91}, already mentioned above, there is up to isomorphisms only a finite number 
of $n$-dimensional lattice polytopes with $l > 0$ interior lattice points, if $n$ and $l$ is fixed.
}
\end{remark}

An indication towards Conjecture \ref{neu} is the following generalization of Corollary \ref{bat-cor}. 
It follows immediately from the previous remark 
and an inequality due to Stanley \cite[Prop.4.1]{Sta91}:
\[1+h^*_1 \leq h^*_{d-1} + h^*_d.\]

\begin{cor}
There is only a finite number of lattice polytopes of fixed degree $d$ and with fixed $h^*_{d-1}$ and $h^*_d$ 
up to isomorphisms and lattice pyramid constructions.
In particular, the volume of any lattice polytope of degree $d$ is bounded by a function depending only on 
$d$, $h^*_{d-1}$ and $h^*_d$.
\end{cor}

\bigskip

The paper is organized in three sections:

In the second section we deal with lattice simplices of degree $d$, showing that they are lattice pyramids over lower-dimensional lattice simplices, if 
their dimension is larger than $4d-2$. Based on this result we prove in the third section Theorem \ref{main-theo} and Proposition \ref{vol-prop}.

\subsection*{Acknowledgments:}
The author would like to thank Christian Haase and Andreas Paffenholz of the Research Group Lattice Polytopes at the Freie Universit\"at Berlin 
for discussions and joint work on this subject.

\section{Lattice simplices with fixed degree}

In this section we prove Theorem \ref{main-theo} for $c=0$:

\begin{thm}
Any lattice simplex of degree $\leq d$ and dimension $n \geq 4d-1$ 
is a lattice pyramid over an $(n-1)$-dimensional lattice simplex.
\label{simplex}
\end{thm}
The bound $4d-1$ is sharp for $d \leq 1$, see \cite{BN07}.\\

Through the whole section let $M = \Z^{n+1}$, and $P = \conv(v_0, \ldots, v_n)$ be an $n$-dimensional lattice simplex of degree $d$, embedded in $\MR = \R^{n+1}$ 
on the affine hyperplane $\R^n \times \{1\}$. We define the half-open {\em parallelepiped}
\[\Pi(P) := \left\{\sum_{i=0}^n \lambda_i v_i \,:\, \lambda_i \in [0,1[ \right\}.\]
Moreover, for $x = \sum_{i=0}^n \lambda_i v_i \in \Pi(P) \cap M$ we define its {\em support}
\[\supp(x) := \{i \in \{0, \ldots, n\} \,:\, \lambda_i \not= 0\},\]
and its {\em height} as the last coordinate of $x$
\[\height(x) := \sum_{i=0}^n \lambda_i \in \N.\]
It is well-known \cite[Cor.3.11]{BR06} that $h^*_i$ equals the number of lattice points in $\Pi(P)$ of height $i$. From this observation, we derive 
the following result:

\begin{lemma}
Let $m \in \Pi(P) \cap M$. Then $\card{\supp(m)} \leq 2 d$.
\label{supp}
\end{lemma}

\begin{proof}

Let $m = \sum_{i=0}^s \lambda_i v_i$ with $\lambda_i \not= 0$ for $i = 0, \ldots, s$. 
We define $P' := \conv(v_0, \ldots, v_s)$. Then $m$ is a lattice point in the relative interior of $\height(m) \cdot P'$, 
so $s+1-\deg(P')=\codeg(P') \leq \height(m)$, hence $s+1 \leq \height(m) + \deg(P')$. 
Since $m \in \Pi(P) \cap M$, we have $\height(m) \leq d$, and by monotonicity $\deg(P') \leq d$. Therefore $\card{\supp(m)} = s+1 \leq 2 d$.
\end{proof}

Let us define the {\em support} of $P$ as
\[\supp(P) := \bigcup_{m \in \Pi(P) \cap M} \supp(m) \subseteq \{0, \ldots, n\}.\]

The relation of this notion to lattice pyramids is straightforward:

\begin{lemma}
Let $i \in \{0, \ldots, n\}$. Then $P$ is a lattice pyramid with apex $v_i$ if and only if $i \not\in \supp(P)$.
\label{pyramid}
\end{lemma}

\begin{proof}

Let $P' := \conv(v_j \,:\, j = 0, \ldots, n,\, j\not=i)$. Then $P$ is a lattice pyramid over $P'$ if and only if $\Vol(P) = \Vol(P')$. 
Now, the statement follows from $\Vol(P) = \card{\Pi(P) \cap M} \geq \card{\Pi(P') \cap M} = \Vol(P')$.
\end{proof}

\medskip

Now, we can give the proof of Theorem \ref{simplex}:

\begin{proof}[Proof of Theorem \ref{simplex}]

By Lemma \ref{pyramid} it is enough to 
show
\[\card{\supp(P)} \leq 4d-1.\] 

Let $m_0 \in \Pi(P) \cap M$ with $I_0 := \supp(m_0)$ maximal. 
Now, we choose successively in a "greedy" manner lattice points $m_0, m_1, \ldots, m_k \in \Pi(P) \cap M$ 
such that $\card{I_k}$ is maximal, where 
\[I_k := \supp(m_k) \backslash \left(\bigcup_{j=0}^{k-1} \supp(m_j)\right).\]

\textbf{Claim:} For $k \in \N$ we have $\card{I_k} \leq \frac{2d}{2^k}$.\\

Assume that the claim were already proven. Then, since $\card{\Pi(P) \cap M}$\\$= \Vol(P)$ is finite, 
the construction yields
\[\card{\supp(P)} = \left|\bigcup_{k=0}^{\card{\Pi(P) \cap M}} \supp(m_k)\right| < \sum_{k=0}^\infty \frac{2d}{2^k} = 4 d.\]
This proves the theorem. It remains to show the claim:\\

The claim holds for $k=0$ by Lemma \ref{supp}. 
Let it be true for $k-1\in \N$. We set $J_k := I_{k-1} \cap \supp(m_k)$. This implies
\[J_k \sqcup I_k \subseteq \supp(m_k) \backslash \left(\bigcup_{j=0}^{k-2} \supp(m_j)\right).\]
Hence, by the choice of $m_{k-1}$ with $I_{k-1}$ maximal we get
\begin{equation}
\card{J_k} + \card{I_k} \leq \card{I_{k-1}}.
\label{eqq1}
\end{equation}
On the other hand, let $m_{k-1} = \sum_{i=0}^n \lambda_i v_i$ and 
$m_k = \sum_{i=0}^n \mu_i v_i$. Now, we translate $m_{k-1}+m_k$ into $\Pi(P)$:
\[m := \sum_{i=0}^n \{\lambda_i + \mu_i\} v_i \in \Pi(P) \cap M,\]
where $\{\gamma\} \in [0,1[$ denotes the fractional part of $\gamma \in \R$. 
By construction, $\mu_i=0$ and $\{\lambda_i + \mu_i\} = \lambda_i > 0$ for $i \in I_{k-1}\backslash J_k$, as well as 
$\lambda_i=0$ and $\{\lambda_i + \mu_i\} = \mu_i > 0$ for $i \in I_k$. This implies 
\[(I_{k-1}\backslash J_k) \sqcup I_k \subseteq \supp(m) \backslash \left(\bigcup_{j=0}^{k-2} \supp(m_j)\right).\]
Again, by the maximality of $\card{I_{k-1}}$ we get
\begin{equation}
\card{I_{k-1}} - \card{J_k} + \card{I_k} \leq \card{I_{k-1}}.
\label{eqq2}
\end{equation}
Combining equations (\ref{eqq1}) and (\ref{eqq2}) yields
\[\card{I_k} \leq \card{J_k} \leq \card{I_{k-1}} - \card{I_k}.\]
Hence, $\card{I_k} \leq \card{I_{k-1}}/2 \leq \frac{2d}{2^k}$ by induction hypothesis. 
This proves the claim.

\end{proof}

\section{Proof of Theorem \ref{main-theo} and Proposition \ref{vol-prop}}

Throughout, let $P \subseteq \MR$ be a lattice polytope of dimension $n$ and degree $\leq d$. 
The proofs here are based on induction. For the induction step we need the notion of a circuit:

\begin{definition}{\rm 
An affinely dependent subset $\calC \subseteq \V(P)$ is called {\em circuit} in $P$, 
if any proper subset of $\calC$ is affinely independent.
}
\end{definition}

The importance of this notion lies in the fact that $P$ is combinatorially a pyramid with apex $v \in \V(P)$ if and only 
if $v$ is not contained in any circuit in $P$.

The following observation \cite[Lemma 2.1]{HNP07} is joint work with Christian Haase and Andreas Paffenholz:

\begin{lemma}
Any circuit in $P$ consists of $\leq 2d+2$ elements.
\label{circuit}
\end{lemma}

\begin{proof}

We may assume as in the previous section that $P$ is embedded in $\R^{n+1}$ on the 
affine hyperplane with last coordinate $1$. In this case, there is a linear relation 
\[\sum_{v \in \calC_1} z_v v = \sum_{w \in \calC_2} z_w w\]
for $\calC = \calC_1 \sqcup \calC_2$ and $z_v, z_w \in \N_{>0}$. Let $Q := \conv(\calC)$. 
The dimension of $Q$ equals $\card{\calC_1} + \card{\calC_2} -2$. 
We observe that $\sum_{v \in \calC_1} v$ is a lattice point in the relative interior of $\card{\calC_1} \cdot Q$. 
Thus, $\codeg(Q) \leq \card{\calC_1}$, so by monotonicity 
$d \geq \deg(Q) = \dim(Q)+1-\codeg(Q) \geq \card{\calC_2} -1$. Hence $\card{\calC_2} \leq d+1$. Symmetrically, $\card{\calC_1} \leq d+1$. 
This proves the statement.
\end{proof}

Using this lemma we can prove Theorem \ref{main-theo}:

\begin{proof}[Proof of Theorem \ref{main-theo}]

First, let us define $n(c,d) := c (2d+1) + 4d-1$. 
Now, we prove by induction on $c \geq 0$ that any $n$-dimensional lattice polytope $P \subseteq \MR$ 
having $\leq c+n+1$ vertices and degree $\leq d$ is a lattice pyramid over a lattice polytope of dimension $< n(c,d)$.

So, let $P$ be given in this way, and $n \geq n(c,d)$.

If $c=0$, then $\card{\V(P)} = n+1$, so $P$ is a simplex, and the statement follows from 
Theorem \ref{simplex}, since $n(0,d) = 4d-1$.

Let $c \geq 1$. Since $P$ is not a simplex, there is a vertex $v \in \V(P)$ such that $Q := \conv(\V(P)\backslash\{v\})$ is an $n$-dimensional 
lattice polytope. Since $(\card{\V(Q)}-n-1) < (\card{\V(P)}-n-1) \leq c$, the induction hypothesis yields that $Q$ is a 
lattice pyramid over a lattice polytope $B$ with $\dim(B) < n(c-1,d)$. 

Now, since $\dim(Q) = \dim(P)$, there is a circuit in $P$ containing vertices $v,w_1, \ldots, w_l$, where 
$w_j \in \V(Q)$ (for $j = 1, \ldots, l$), and $l \leq 2d+1$ by Lemma \ref{circuit}. In particular, 
$v \in \aff(w_1, \ldots, w_l)$. We set $D := \conv(B, w_1, \ldots, w_l) \subseteq Q$. 
Hence, $Q$ is a lattice pyramid over the lattice polytope $D$, whose dimension satisfies
\[\dim(D) \leq \dim(B) + l \leq (n(c-1,d)-1) + (2d+1) = n(c,d)-1.\]
 Since 
$\aff(D) = \aff(D,v)$, also $P$ is a lattice pyramid over the lattice polytope $\conv(D,v)$ of dimension $\dim(D) < n(c,d)$. 

\end{proof}

The proof of Proposition \ref{vol-prop} is analogous:

\begin{proof}[Proof of Proposition \ref{vol-prop}]

The proof is by induction on $V \geq 1$. Let $P \subseteq \MR$ be a lattice polytope 
having volume $V$, degree $d$, and $\dim(P) = n \geq (V-1) (2d+1)$. 

If $V=1$, the statement is trivial. So, let $V \geq 2$. 

First, let $P$ be a lattice simplex. If $V \geq 3$, then $n \geq 4d+2$, so the statement follows from Theorem \ref{simplex}. 
If $V = 2$, then there exists in the notation of the previous section precisely one lattice point $0 \not= m \in \Pi(P) \cap M$. 
Hence, $\card{\supp(P)} = \card{\supp(m)} \leq 2d$ by Lemma \ref{supp}, so $P$ is a 
lattice pyramid over an $(n-1)$-dimensional lattice simplex by Lemma \ref{pyramid}, since $n \geq 2d+1$. 

Therefore, we can assume that $P$ is not a simplex. Now, the remaining induction step proceeds precisely 
as in the proof of Theorem \ref{main-theo}. 
\end{proof}

\end{document}